\documentclass[runningheads]{llncs}
\usepackage{graphicx}
\usepackage[english]{babel} 
\usepackage[T1]{fontenc} 
\usepackage{ae} 
\usepackage{amssymb} 
\usepackage{graphicx} 
\usepackage[utf8]{inputenc} 
\usepackage{graphicx} 
\usepackage{amsmath,amsfonts,amssymb,dsfont,mathtools,blindtext} 
\usepackage{changepage}
\usepackage{longtable}

\usepackage{subfigure}
\usepackage{multirow}
\usepackage{adjustbox}
\usepackage{cite}

\begin{document}
\title{Planning routes in road freight minimizing logistical costs and accident risks}
%
%
\author{Gabriel Adam Bilato\inst{1}\orcidID{0000-0002-5427-3958} \and
Cleber Damião Rocco\inst{2}\orcidID{0000-0002-7988-6136} \and
Anibal Tavares de Azevedo\inst{3}\orcidID{2222--3333-4444-5555}}
\authorrunning{G. Bilato et al.}
%
\institute{Universidade Estadual de Campinas, Limeira SP, Brazil \email{g168460@g.unicamp.br}\\ \and
Universidade Estadual de Campinas, Limeira SP, Brazil
\email{cdrocco@unicamp.br}\\
\and
Universidade Estadual de Campinas, Limeira SP, Brazil\\
\email{atanibal@unicamp.br}}
\maketitle              
\begin{abstract}
The Vehicle Routing Problem (VRP) has been widely studied throughout its history as a way of optimizing routes by minimizing distances, and the issue of risk in VRP has been received less attention, which is essential to increase transport safety, to reduce accident costs and to improve delivery reliability. In this way, this paper aims to support decision makers to plan routes for a road freight company considering both, i.e. logistics cost and safety.
An analytical approach based on statistics was developed in which official government data of accidents occurrences and data from cargo insurance companies were used to estimate the risk cost of routes using the Monte Carlo simulation.
The Capacitated Vehicle Routing Problem (CVRP) was employed and logistics and risk costs were minimized by varying a specific safety level coefficient and the model generated solutions with safer routes, reducing risk cost by up to approximately 18\%. The accidents probabilities and risk costs of each route represented values as expected and the main contributions of this paper is the applicability of the approach to support decision markers to choose routes considering safety and logistics costs, and to be a simple and adaptable methodology for any VRP model. In addition \textit{Knime Analytics Platform} was also used to estimate the accidents probabilities and to simplify data exploration, analysis, visualization and interpretation.

\end{abstract}
\section{Introduction}
In an increasingly globalized and connected world in which the population and urbanization grows over the years, the pressure for more efficient and sustainable road freight increases due to the greater demand for different types of products. Road freight plays a fundamental role in global logistics due to the flexibility of the infrastructure that this segment of transport offers, allowing a door-to-door service, which is not possible with other modes \cite{Ergstron_2016}.

On the other hand, the reliability in road freight depends on some factors that can be the cause of many accidents, interrupting the supply chain and generating significant losses. According to the World Health Organization (WHO) data, accidents costs on average around 3\% of the Gross Domestic Product (GDP) of a country and led to the death of 1.35 million people, being the eighth leading cause of death all around the world.

Therefore, reducing the accident rate of road freight is an important factor from a strategic and sustainable point of view, which seeks logistical and economical development, because when an accident occurs losses can assume large financial proportions in addition to impact in others different dimensions. From this point, the need for improvement in the reduction of losses in road transport stands out \cite{Ergstron_2016}.

When a route is chosen for the vehicle that will leave a depot to deliver to other cities, not only the logistics costs or distance should be considered, but also the measurement of the risks linked to that route. Many studies consider only the first factor using mathematical models like VRP as a way of solving it. Risk measurement is addressed in cash-in-transit problems also through the VRP, which basically is a mathematical model that aims to minimize the distances traveled with the restriction that the value of the risk of robberies of heavy trucks during the transport of money is limited by a risk threshold \cite{talarico_et_al_2015}.

The route safety is also discussed in studies of hazardous materials transportation such as fuels, flammable materials, gases and others, aiming to reduce social and environmental impacts and increase transport safety by minimizing the risk factor \cite{holeczek_2021}.

VRP has been widely studied throughout its sixty-year history, but still few studies have taken into account the issue of risk. Recently, researches have emerged in the areas of cash-in-transit and hazardous materials, but they did not consider a statistical approach to risk into the models.


Therefore, verifying the applications of the topic risks in VRP and the need to reduce losses due to accidents, this paper aims to introduce an analytical approach for road freight through optimization model and statistical analysis to support the decision making for the choice of routes based on logistic cost and safety.

In addition to the contributions of this analytical approach are its simplicity of application to a real problem of a Brazilian road freight company that daily needs to choose the best routes that consider safety and costs, as well as their adaptability to other VRP models. Finally, Knime Analytics Platform was helpful to simplify data exploration, analysis, visualization and interpretation and to estimate the accidents probabilities.

This paper is organized in five sessions in which the first one introduces the importance of studying risk in VRP. From the second session, a review and discussion of the literature which addressed risks in VRP is carried out. In an analytical approach, the methodology proposed for the study is shown, as well as the calculations used. Finally, the results are presented in experimental studies and the conclusion is drawn from this.

\section{Literature review}
\subsection{Risk in VRP}

Vehicle Routing Problems (VRP) have been extensively studied throughout their history to support real-life applications. Risk and safety in VRP has been received more attention in applications for the transport of hazardous materials, whose risk is an accident causing socio-environmental damage, and cash-in-transit which are related to cargo theft \cite{talarico_et_al_2017}.

\cite{erkut_Ingolfsson_2005} cited eight risk models that were developed for the optimization of hazardous materials transportation, whose three main ones are represented by Table \ref{tab:table_1}. For the first model, the risk ($IP$) is calculated through the probability of the undesirable occurrence event of each route $i$, while the second ($PE$) considers only the number of people exposed to risk. In the traditional model, the risk ($TR$) is a product between the probability of the undesirable event and its measure of consequence.

To investigate the behavior of these three risk models in VRP, \cite{holeczek_2021} used bi-objective functions that minimize both distance and accident risk. The results shows that the Traditional Risk generates the best reduction in total risk but it shows the worst deviation from the minimum distance when compared to the other two models. The Accident Probability offers the best trade-off with an economical goal and it is most appropriate for problems where the consequences are uncertain. As for the Population Exposed, the data are more easily acquired and the results are evaluated more intuitively for a decision maker, but it could only be applied to problems that consider urban areas because for an environment such as rural areas other factors must be considered.

\begin{table}[]
    \centering
    \begin{tabular}{l l l}
    \textbf{Model}            & \textbf{Equation} &  \\\hline
    Accident Probability  &  $IP(r) = \sum_{i\in r} p_i$ & $p_i =$ accident probability \\\hline
    Population Exposure         &  $PE(r) = \sum_{i\in r} D_i$ & $D_i =$ population exposure \\\hline
    Tradicional Risk                &  $TR(r) = \sum_{i\in r} p_i * C_i$ & $p_i =$ accident probability \\
     & & $C_i =$ measure of the consequence \\\hline
    \end{tabular}
    \caption{Models for risk assessment. Adapted from \cite{erkut_Ingolfsson_2005}.}
    \label{tab:table_1}
\end{table}

The elements of each model such as accident probability, consequences, number of people exposed and others, are defined according to the case is being studied. \cite{androutsopoulos_et_al_2012,carrese_et_al_2019} , for example, studied the risks in hazardous materials transport and they considered that the undesirable event would be the accident whose consequence is related to the number of people exposed to risk. On the other hand, \cite{talarico_et_al_2015} applied the risk in \textit{cash-in-transit} routing problem that they considered the probability of a robbery is proportional to the distance of the route and, as a consequence, the amount of cash into the truck.

In addition to the elements, the models also vary according to the risk approach used in VRP. Thus, Table \ref{tab:table_2} was built to show the case studied, models and risk elements that were adopted for the calculation. 

\cite{du2017multi,Pradhananga2014,Wang2018} applied the concept of the traditional model, as proposed by \cite{erkut_Ingolfsson_2005} in which the probability of an accident was used as an undesirable event and the exposed population as a consequence. \cite{androutsopoulos_et_al_2012,holeczek_2021} also applied the traditional definition, however it was considered as \textit{load dependent}, that is, the amount of load factor is added to the model and varies as deliveries are made. \cite{carrese_et_al_2019} is also based on the traditional model, but two other factors that interfere in the driver's attention are added to the objective function, the Altimetric Index and the Planimetric Index. The first considers the elevations along the route while the second is introduced to take into account geometrical constraints related to the road radius.

\cite{Bula_et_al_2016,Bula2019} calculated the risk in a different way than what has been discussed so far. In this case, in addition to the accidents probabilities, it also considered the release probability as a result of the accident, which the consequences are the number of people exposed, and the truck type and the load amount have a relevant impact on the risk.

As a derivation of hazardous materials, models for \textit{cash-in-transit} arise. \cite{talarico_et_al_2015,talarico_et_al_2017} proposed the traditional method to calculate the risk of robbery. As already explained, they considered the consequence as equivalent to the amount of value in transport. \cite{Ghannadpour2020} also assumed the consequence as the same way and the distance is proportional to the risk of theft, but added to the model a factor about the frequency of passing through the same route and the ambushing probabilities to the vehicle and its success.

Table \ref{tab:table_2} shows that the studies are concentrated in only two areas: hazardous materials and cash-in-transit. Talking about other segments of transport, there are few studies in the literature taking into account risks in VRP and it is important because there are also significant consequences or losses in case an accident happens, especially financial when the values of shipped goods are high.

The VRP model with risks vary according to the objectives of each study. Table \ref{tab:table_3} summarizes the problem characteristics with the nomenclature presented by \cite{Braekers2016}. The symbol ``x'' is used when the paper considers the characteristics and Table \ref{tab:my_label} represents the nomenclatures.

\cite{Bula_et_al_2016,talarico_et_al_2015} used a mono-objective function which firstly minimized the distances and than the risk. \cite{androutsopoulos_et_al_2012,Bula2019,Pradhananga2014,Wang2018} suggest bi-objective functions that analyze both logistical costs or distances and route risks. Multi-objectives are presented by \cite{carrese_et_al_2019,Ghannadpour2020,Zheng2010} which the first one, in addition to minimize the distance, weighted the risk by the traditional method (probability of accident and a consequence) and the number of people exposed.

Generally, the risk factor is analyzed as an objective to be minimized, however \cite{talarico_et_al_2015} considers it as a constraint which the risk value is limited by a risk threshold and classified as \textit{Risk constrained Cash-in-Transit Vehicle Routing Problem (RCTVRP)}. \cite{Wang2018} restricts the condition that no vehicles of the same fleet travel in echelon because when there are two or more vehicles using the same route as the same time the consequences are considered to be greater if occurs an accident between them.

\begin{longtable}{p{.30\textwidth} p{.30\textwidth} p{.30\textwidth}}
    \textbf{Authors} & \textbf{Case studied} & \textbf{Risk model}  \\\hline

    \cite{Zheng2010} & Hazardous Materials & accident probability * consequence + population exposure \\\hline
    
    \cite{androutsopoulos_et_al_2012,holeczek_2021} & Hazardous Materials & accident probability * population exposure * load amount \\\hline
    
    \cite{du2017multi,Pradhananga2014,Wang2018} & Hazardous Materials & accident probability * population exposure \\\hline
    
    \cite{talarico_et_al_2015} & Cash-in-Transit & route length * load amount \\\hline
    
    \cite{Bula_et_al_2016,Bula2019} & Hazardous Materials & accident probability * release probability * route length * load amount * truck type * population exposure \\\hline
    
    \cite{carrese_et_al_2019} & Hazardous Materials & accident probability * population exposure + altimetric index + planimetric index \\\hline
    
    \cite{Ghannadpour2020} & Cash-in-Transit & ambush probability * theft success probability * route length * load amount * frequency of repeated use of a route \\\hline

    \caption{Risk models proposed in literature.}
    \label{tab:table_2}
\end{longtable}

\begin{table}[]
    \centering
    \begin{tabular}{c c}
      \multicolumn{2}{l}{}\\   
      \hline
      \textbf{Categories}  &  \textbf{Sub-categories} \\\hline
      \multirow{4}{*}{Objective Function (3.10)} & Travel time dependent (3.10.1) \\
                                & Distance dependent (3.10.2)\\
                                & Implied hazard/risk related (3.10.5)\\
                                & Others (3.10.6)\\\hline
      Data used (5.1)           & Real-world data (5.1.1) \\\hline
    \end{tabular}
    \caption{Nomenclature proposed by \cite{Braekers2016}.}
    \label{tab:my_label}
\end{table}

\begin{table}[]
    \centering
    \begin{tabular}{l c c c c c c}
    \multicolumn{5}{l}{\small{* Nomenclature presented by \cite{Braekers2016}.}}\\
    \hline
      \textbf{Authors}  & \textbf{VRP} & \multicolumn{4}{c}{\textbf{3.10*}} & \textbf{5.1*}\\
    \hline
      \textbf{} & \textbf{} & \textbf{3.10.1*} & \textbf{3.10.2*} & \textbf{3.10.5*} & \textbf{3.10.6*} & \textbf{5.1.1*} \\\hline
      
      \cite{androutsopoulos_et_al_2012} & VRPTW & x & & x &  \\\hline
      
      \cite{Bula_et_al_2016} & HVRP & & & x & \\\hline
      
      \cite{Bula2019} & HVRP & & x & x & \\\hline
      
      \cite{carrese_et_al_2019} & VRPTW & x & & x & x & x \\\hline
      
      \cite{du2017multi} & MRVRP & & & x & & x \\\hline
      
      \cite{Ghannadpour2020} & VRPTW & & x & x & \\\hline
      
      \cite{holeczek_2021} & CVRP & & x & x & \\\hline
      
      \cite{Pradhananga2014} & VRPTW & x & & x & & x \\\hline
      
      \cite{talarico_et_al_2015} & RCTVRP & & x & & \\\hline
      
      \cite{Wang2018} & VRPTW & & x & x & \\\hline
      
      \cite{Zheng2010} & CVRP &  & x & x & x \\\hline
    \end{tabular}
    \caption{Description and characteristics of VRP models.}
    \label{tab:table_3}
\end{table}

For the measurement of risk, few studies is addressed on how the data are explored and, according to \cite{du2017multi}, it is necessary to integrate real historical data of accidents and big-data in the formulation of models for the transport of hazardous materials, however \cite{talarico_et_al_2015} mentions that there is little data available and \cite{androutsopoulos_et_al_2012} does not explore risk measurement due to its complexity and also states that future studies should deal with this issue.

\cite{carrese_et_al_2019} calculated the accident probability from data obtained by the mobility agency in Rome, quantified population density through a \textit{census data} and measured the infrastructure through the \textit{Google Application Programming Interface (API)}.

\cite{Pradhananga2014} estimated accident rates using data collected from the \textit{Institute for Traffic Accident Research and Data Analysis} (ITARDA) and the \textit{Ministry of Land, Infrastructure, Transport and Tourism} (MLIT), both from Japan. \cite{Ghannadpour2020} estimated the probability of a robber attack using game theory, while to calculate the success probability it used the multi-criteria decision making.

For a future study, \cite{Pradhananga2014} proposed extensions of the model considering such characteristics for hazardous materials routing problem using real-time traffic information and considering the effects of infrastructural characteristics of the road network for future studies.

\cite{milovanovic_2012} developed a methodology for calculating the risk of accidents in the hazardous materials transport which factors that influence the accident probability as well as their consequences were considered. The factors were measured from indirect interviews with experts obtaining numerical risk results for each route through this analysis. On the other hand, it did not use mathematical models of VRP to optimize routes and it did not consider statistical analysis.

As much as some studies still try to deal with the use of real data in their problems, the statistical view in data analysis is not approached in a deeper way. \cite{Fillbrunn2017} reviewed some the extensions of the free software Knime Analytics Platform that could support the analyzes from a database and that provides the creation of structured workflows. Combined with this, \cite{Ali2021} used the Monte Carlo simulation to assess the losses in economic values for the Pakistani economy in the event of a transport strike, but in the case of this study it would relate the losses to road accidents.

\subsection{Contributions of this paper}

\begin{itemize}
    \item Develop a risk model in VRP that is applied to both the cases studied presented in Table \ref{tab:table_2} and to different cases of the road freight;
    \item Use real-world data and statistical approach to estimate the accident probabilities. As Table \ref{tab:table_3}, few studies take into account real-world data and none of them use a statistical approach.
\end{itemize}

\section{Analytical Approach}
\subsection{Description}

\begin{figure}
    \centering
    \includegraphics[width=1\textwidth]{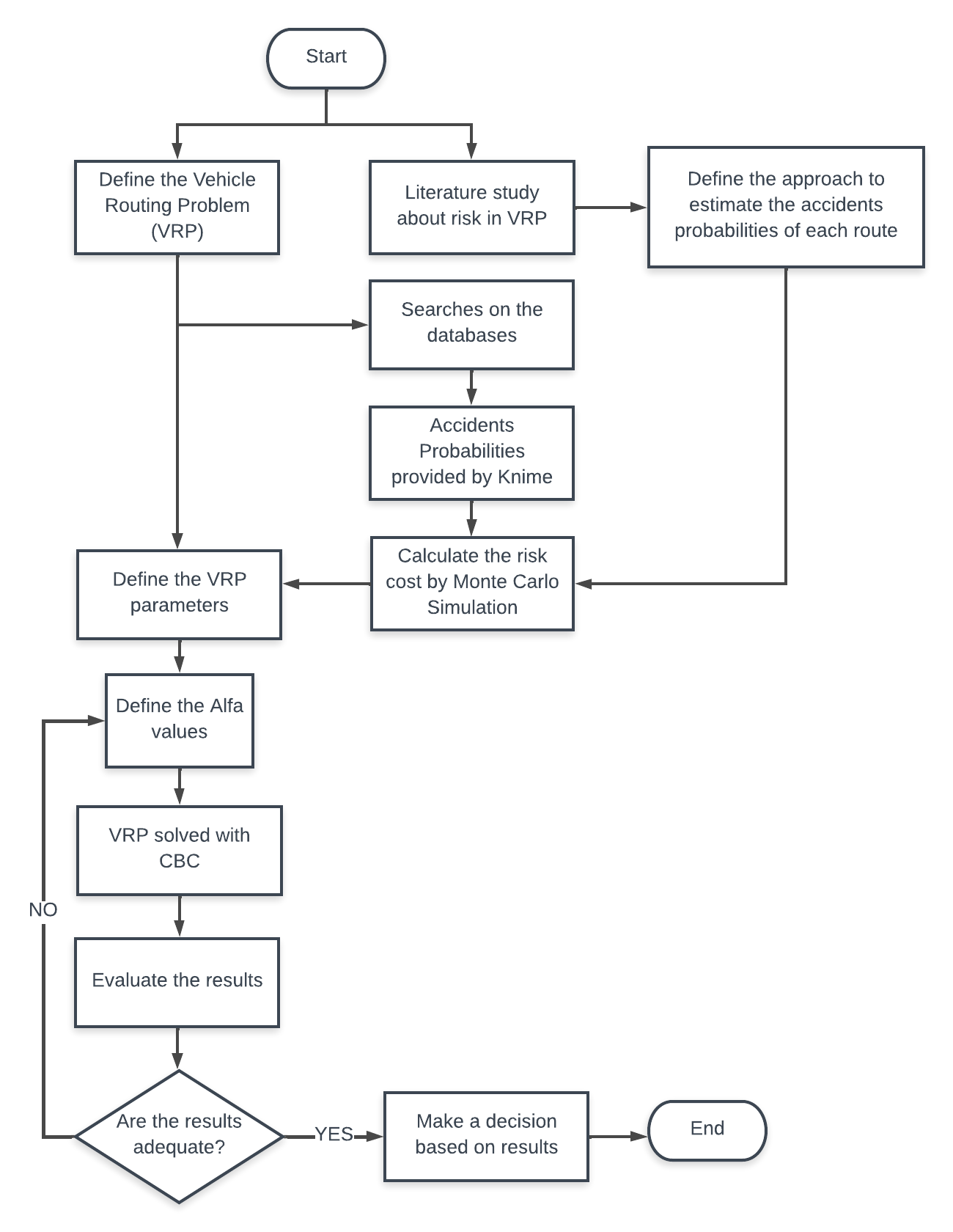}
    \caption{Workflow to employ this paper metodology.}
    \label{figure1}
\end{figure}

This study followed the procedures described in the workflow (Figure \ref{figure1}) and started from the VRP definition that was employed with a set of cities and roads known by the authors. Then, the parameters were defined by the online tool \cite{qualp} to calculate the logistics costs of each arc, which consider fuel expenses, based on the vehicle's consumption and its value per liter, and with tolls, whether there is on the arc.

In parallel, the literature on statistical analysis of risk in road freight was revised and from \cite{milovanovic_2012} and from the data that was found, the approach to calculate the probabilities and costs related to the risk of accident for each arc was developed. All collected and processed data used in this study is available at \cite{github}. 

Several databases of Brazilian roads were consulted on the \textit{web} and it was found that, for the most part, it comes from government agencies whose information is available for public viewing and among them there are: the National Transport Confederation (CNT), the Department of Roads and Highways of the State of Sao Paulo (DER- SP) and the National Department of Transport Infrastructure (DNIT). Through a load insurance, data was provided regarding the losses in monetary values that occurred due to accidents between the periods from January 2018 to March 2021.

The data were processed using the \textit{Knime Analytics Platform} tool, which generated the accident probabilities for each arc. Then, from the probability results, the Monte Carlo simulation was programmed to obtain the costs related to the risks of the arcs.


The objective function (Equation \ref{equation1}) minimizes logistics and risk-related costs, and it was weighted by the $\alpha$ parameter whose values vary between 0 and 1 and which consists of assigning weight to safety during the optimization.

After defining the objective function and the parameters of the VRP, the mathematical model was implemented in \textit{Python} and the results were generated by CBC solver. The decision maker could evaluate the routes and the values of each objective function component, logistics cost and risk cost, and for the total that is the sum of them.

The security level $\alpha$ must be adjusted by the decision maker to assess whether the results are suitable to the case. Depending on the cargo type or its value, $\alpha$ may have to be adjusted in order to have a greater or lesser level of safety prioritization of the routes.

Therefore, this evaluation is necessary for the decision maker to consider logistics costs and security, if the result is not adequate, a new change of $\alpha$ must be carried out until the solution fits to the best choice of the decision maker.

\subsection{The VRP model}

\begin{figure}[h]
    \centering
    \includegraphics[width=0.8\textwidth]{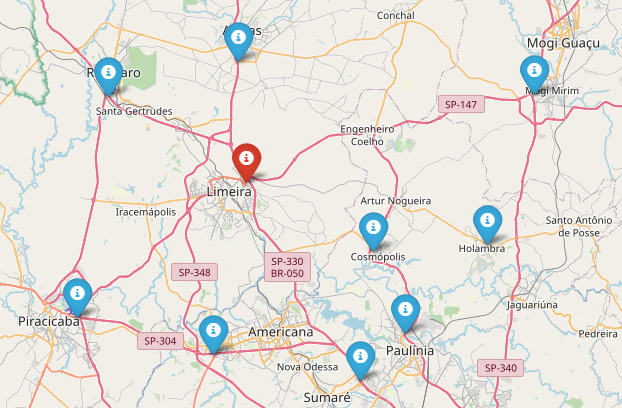}
    \caption{Depot in red and delivery points in blue.}
    \label{figure2}
\end{figure}

The capacitated VRP (CVRP) is employed in a real problem of a transportation company and the arcs are defined through the indirect graph $G = (N,E)$. The problem presents a single depot that is located in Limeira/SP. Two vertices are created to represent exit and arrival in the depot, respectively ($\{0,n+1\}$). Set $C= \{1,...,n\}$ represents the delivery points in nine cities ($n=9$). The total number of vertices is represented by $N = C \cup \{0,n+1\}$. Set ($K = \{k_1,k_2,k_3\}$) is the vehicles with equal capacity $q$. The map indicated in Figure \ref{figure2} illustrates both the depot and all the delivery points that the carrier must carry out.

The expressions of the mathematical model are presented next. Firstly it is defined the sets.

$N = C \cup \ \{0,n+1\}, C= \{1,...,n\}, E = \{(i,j) : i,j \in N, i \neq j, i \neq n+1, j \neq 0\}$
\\
\begin{equation}
   min z = (1 - \alpha)\sum_{k\in K}\sum_{(ij)\in E}c_{ij}X_{ijk} + \alpha\sum_{k\in K}\sum_{(ij)\in E}r_{ij}X_{ijk}
    \label{equation1}
\end{equation}

\begin{equation}
   \sum_{k\in K}\sum_{j\in E}X_{ijk} = 1 , \forall i \in C
    \label{equation2}
\end{equation}

\begin{equation}
   \sum_{i\in E}d_i\sum_{j\in E}X_{ijk} \leq q , \forall k \in K
    \label{equation3}
\end{equation}

\begin{equation}
   \sum_{j\in E}X_{0jk} = 1 , \forall k \in K
    \label{equation4}
\end{equation}

\begin{equation}
   \sum_{i\in E}X_{ihk} - \sum_{j\in E}X_{hjk} = 0 , \forall h \in C, k \in K
    \label{equation5}
\end{equation}

\begin{equation}
   \sum_{i\in E}X_{i,n+1,k} = 1 , \forall k \in K
    \label{equation6}
\end{equation}

\begin{equation}
   u_{ik}-(n+1)X_{ijk} \leq u_{jk} - n, \forall (i,j) \in E, k \in K
    \label{equation7}
\end{equation}

\begin{equation}
   X_{ijk} \in \{0,1\}, \forall (i,j) \in E, k \in K
    \label{equation9}
\end{equation}

Parameters: $c_{ij}$: logistics cost; $r_{ij}$: Risk cost; $\alpha$: security level; $d_i$: demand of node $i$; $q$: vehicle capacity; $X_{ijk}$: binary variable.

The Equation (\ref{equation1}) is the objective function that minimizes logistics $c_{ij}$ and risks $r_{ij}$ costs. The Equation (\ref{equation2}) describes each arc must receive only one vehicle. The vehicle capacity constraint is represented by the Expressions (\ref{equation3}). Equations (\ref{equation4}) and (\ref{equation6}) show every vehicle must leave and arrive at the depot, respectively. The flow restriction for each node is represented by the Equation (\ref{equation5}). Finally, the sub-tour elimination is indicated by the Constraints (\ref{equation7}) and the binary variables are defined by Expression (\ref{equation9}).

Expressions (\ref{equation8}) ensure a strong formulation to eliminate sub-tour because of its LP relaxation but the number of these constrains could be too high ($2^n$) to the model. Because of this, it was used \textit{cut callbacks} and \textit{lazy constraints} to insert only constrains which violated the model and it is not necessary to use $2 \leq |S| \leq n$ in the model \cite{HaroldoG.Santos2019}.

\begin{equation}
   \sum_{ij \in S}{X_{ijk}} \leq |S| -1, \forall S \subseteq C, k \in K
    \label{equation8}
\end{equation}

Where: $S$ is a set of sub-tour.


The \textit{cut callbacks} are used only to improve the LP relaxation but not to define feasible solutions, which need to define by the initial formulation. So, it is included a weak sub-tour elimination constraints presented by Expression (\ref{equation7}) in the initial model and then add Constraints (\ref{equation8}) as cuts.

In this way, \textit{Python MIP} (\cite{HaroldoG.Santos2019}) package was introduced to the model because it is possible to eliminate sub-tour with \textit{cut callbacks} and \textit{lazy constraints} that makes it more efficient. The free solver CBC was used to optimize this problem and it was solved in 44 seconds for all twenty different values of $\alpha$ resulting in 363 variables and 289 constrains for each instance.


\subsection{Calculating the risk cost $r_{ij}$}
\begin{figure}
    \centering
    \includegraphics[width=0.75\textwidth]{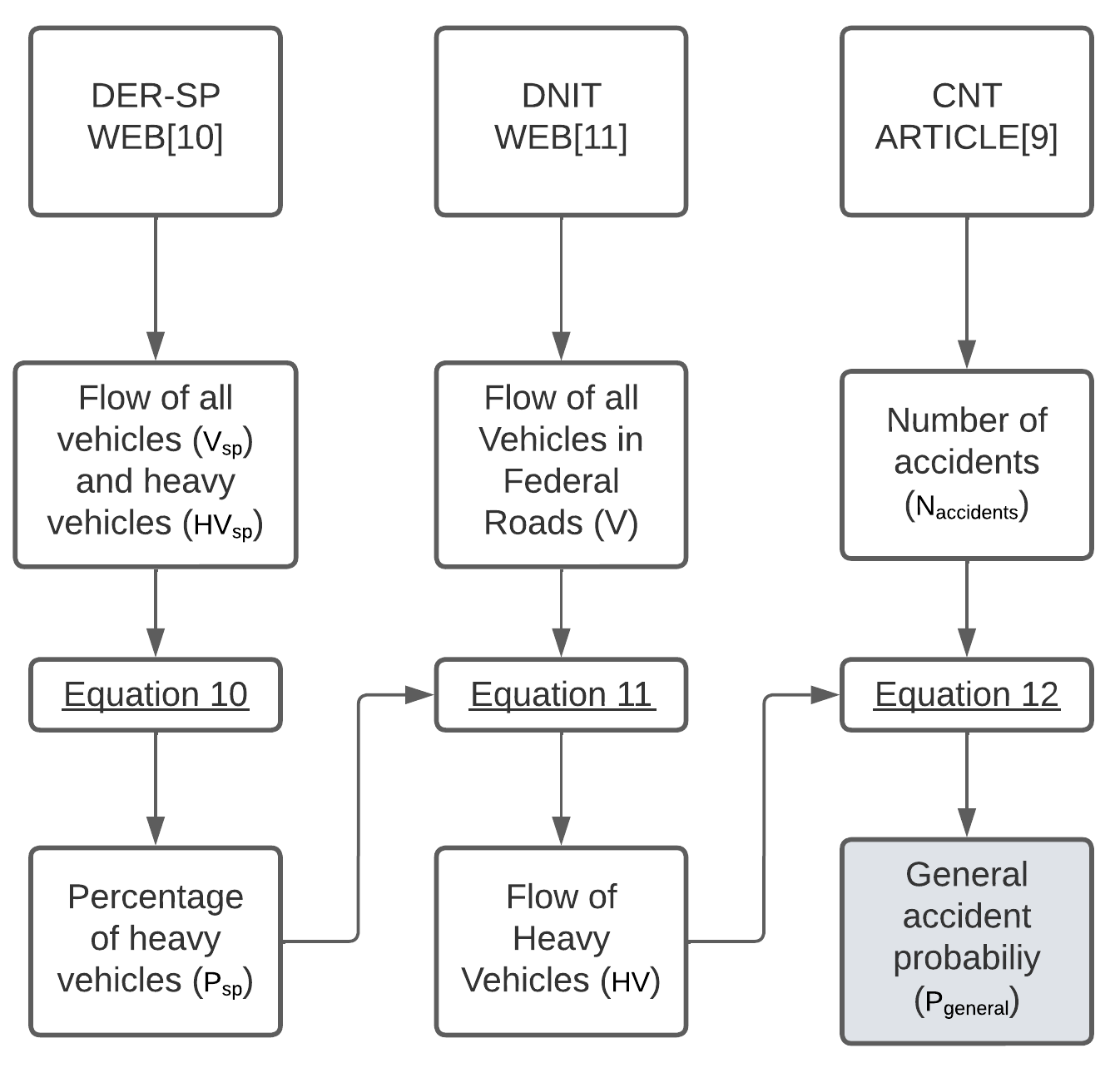}
    \caption{Workflow to calculate $Pgeneral$.}
    \label{figure3}
\end{figure}

\begin{figure}
    \centering
    \includegraphics[width=0.9\textwidth]{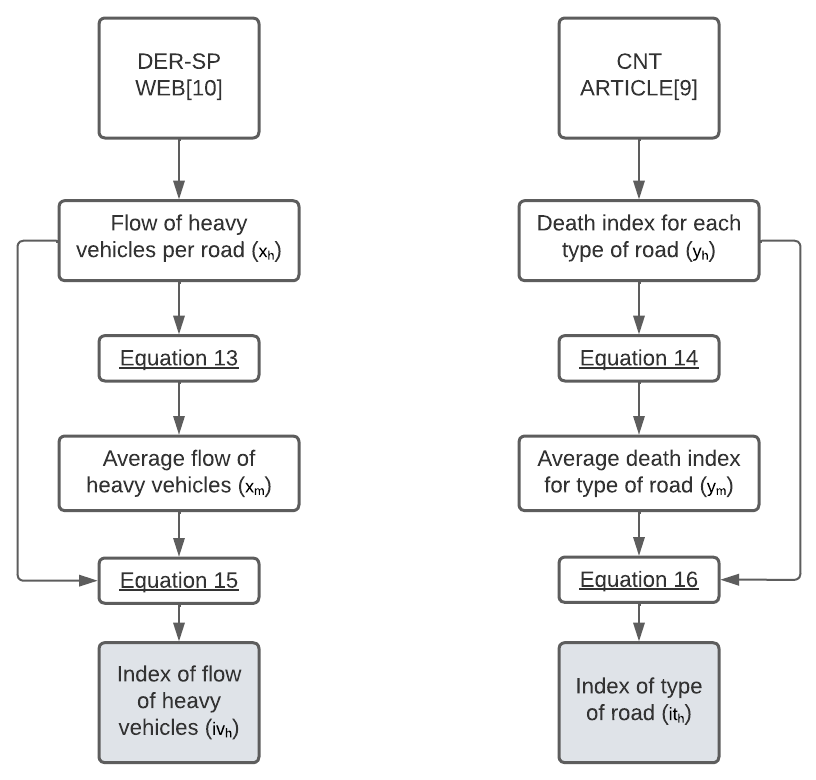}
    \caption{Workflow to calculate the indexes $iv_h$ e $it_h$.}
    \label{figure4}
\end{figure}

\begin{figure}
    \centering
    \includegraphics[width=0.9\textwidth]{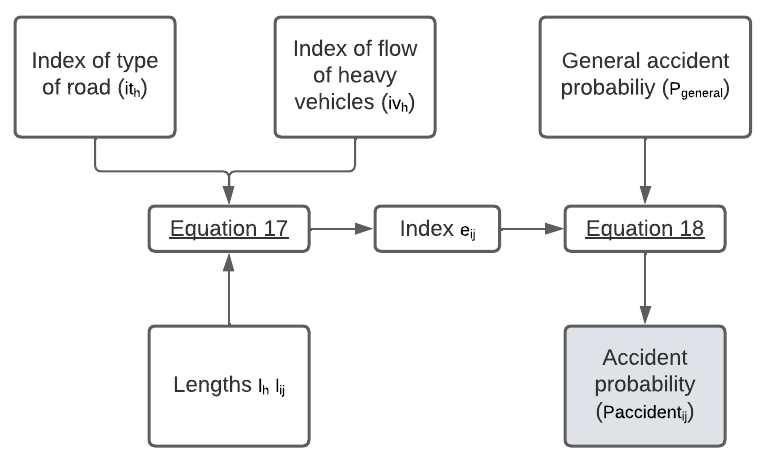}
    \caption{Workflow to calculate $e_{ij}$ and $Paccident_{ij}$.}
    \label{figure5}
\end{figure}

As already mentioned, the accidents probabilities $Paccident_{ij}$ were generated for each arc \textit{ij} by \textit{Knime Analytics Platform}, according to the workflows represented by the Figures \ref{figure3}, \ref{figure4} and \ref{figure5} and Equations (\ref{equation10}), (\ref{equation11}), (\ref{equation12}), (\ref{equation13}), (\ref{equation14}), (\ref{equation15}), (\ref{equation16}), (\ref{equation17}) e (\ref{equation18}).


At this first moment, it was necessary to estimate a general probability ($Pgeneral$) of accidents occurring on any road in Brazil according to the workflow of the Figure \ref{figure3}. Data was collected from the free websites of \cite{dersp,dnit} and from the article \cite{cnt_painelacidentes2018}. \cite{dnit} provides the average volume $V$ of all vehicles that travel daily on federal roads, while \cite{dersp} extracted the percentage $P_{sp}$ of heavy vehicles $HV_{sp}$ over the total $V_{sp}$ that circulates in Sao Paulo State roads that was calculated by Equation (\ref{equation10}) and assumed that the percentage in federal roads as the same. From this point of view, Equation (\ref{equation11}) resulted in the amount of heavy vehicles $HV$ that is used in Equation (\ref{equation12}) with the number of accidents $N_{accidents}$, extracted from \cite{cnt_painelacidentes2018}, and then $Pgeneral$ was found.


\begin{equation}
    P_{sp} = HV_{sp}/V_{sp}
    \label{equation10}
\end{equation}

\begin{equation}
    HV = P_{sp}.V
    \label{equation11}
\end{equation}

\begin{equation}
    Pgeneral = \frac{N_{accidents}}{HV}.100\%
    \label{equation12}
\end{equation}

The factors considered to the calculation of $r_{ij}$ were: the types of roads and the traffic of heavy vehicles. The first is considered because in Brazil there are several types of roads that show different safety levels. The Accident Panel Report prepared by \cite{cnt_painelacidentes2018} breaks it into five categories ( $T=\{t1,...,t5\}$) as presented by Table \ref{tab:table_4}. It is also shown that the number of deaths per hundred accidents was extracted from this report and it is considered in the calculation for the risk of accidents.

\begin{table}[]
    \centering
    \begin{tabular}{c c}
      \textbf{Road type}  & \textbf{Death rate per 100 accidents} \\\hline
      two lanes two way road with central safety lane      &  12.3 \\\hline
      two lanes two way road with central barrier   &  8.5  \\\hline
      two lanes two way road with central line      &  18.0 \\\hline
      Single lane one way road                &  11.9 \\\hline
      Single lane two way road                &  22.3 \\\hline
    \end{tabular}    \caption{Death rate per type of road \cite{cnt_painelacidentes2018}.}
    \label{tab:table_4}
\end{table}

The flow of heavy vehicles was obtained from \cite{dersp} by speed radars installed on the roads and it was considered as directly proportional to the accidents probabilities because two roads with the same characteristics but one have a bigger flow than other so the first has greater chances of an accident occurring.

The indexes $it_h$ and $iv_h$ represents, respectively, the types of road and the flow of vehicles were calculated according to the workflow of the Figure \ref{figure4} and where $h \in H$ which means that road $h$ is on the set $H$ of problems roads. The data was collected from \cite{cnt_painelacidentes2018,dersp}. First, it was necessary to find $\Bar{x}$ and $\Bar{y}$ by the Equations (\ref{equation13}) and (\ref{equation14}) that represent, respectively, the average flow of vehicles $x_h$ running on the twelve roads of the problem $N_h = 12$ and the average death rate per 100 accidents $y_t$ among the five types of roads $N_t = 5$. Then, the indexes $iv_h$ and $it_h$ were calculated by the Equations (\ref{equation15}) and (\ref{equation16}), which basically consists if $iv_h$ or $it_h$ is greater than 1.0, the accident probability on the road $h$ will be greater than the general probability ($P_{general}$), and otherwise if the indices are less than 1.0.

\begin{equation}
    \Bar{x} = \frac{\sum_{h\in H}x_{h}}{N_h}
    \label{equation13}
\end{equation}

\begin{equation}
    \Bar{y} = \frac{\sum_{t\in T}y_{t}}{N_t}
    \label{equation14}
\end{equation}

\begin{equation}
    iv_h = 1 + \frac{x_{h} - \Bar{x}}{\Bar{x}} \ , \ \forall h \in H
    \label{equation15}
\end{equation}

\begin{equation}
    it_h = 1 + \frac{y_{h} - \Bar{y}}{\Bar{y}} \ , \ \forall h \in H
    \label{equation16}
\end{equation}

Where $x_h$ means the flow of vehicles represented by road $h$ and $y_h$ the death rate of the type of road $h$.

In some arcs, the vehicle may pass through more than one road ($h$), e. g., the arc Limeira-Cosmópolis has two roads $h_1$=SP330 and $h_2$=SP133 whose flows and characteristics are different. Thus, it is necessary a index $e_{ij}$ that is weighted by $iv_h$, $it_h$ and $l_h$, and represented by the Equation (\ref{equation17}). Where $l_h$ is the length that the truck travels in each road $h$ of the arc $ij$ and $l_{ij}$ represents the total length of the arc $ij$.

\begin{equation}
    e_{ij} = \frac{\sum_{h\in H}iv_h.it_h.l_h}{l_{ij}} \ , \ \forall (i,j) \in E
    \label{equation17}
\end{equation}

Finally, Equation (\ref{equation18}) describes the accident $Paccident_{ij}$ occurring among $ij$. The workflow of the Figure \ref{figure5} shows how $Paccident_{ij}$ and $e_{ij}$ were found.

\begin{equation} 
    Paccident_{ij} = Pgeneral.e_{ij} \ , \ \forall (i,j) \in E
    \label{equation18}
\end{equation}

After to find $Paccident_{ij}$, it is possible to estimate $r_{ij}$ by Monte Carlo simulation. This was possible through data provided by the load insurance, which collected the percentages of accidents divided into intervals of losses in load values according to Table \ref{tab:table_5}.

\begin{table}[]
    \centering
    \begin{tabular}{c c}
      Range of values  & occurrence \\\hline
      \$ 0.01 to \$ 200,000.00                 &  37.91\% \\\hline
      \$ 200,000.00 to \$ 300,000.00           &  24.17\%  \\\hline
      \$ 300,000.00 to \$ 500,000.00           &  19.91\% \\\hline
      \$ 500,000.00 to \$ 1,000,000.00         &  16.11\% \\\hline
      \$ 1,000,000.00 or more                   &  1.90\% \\\hline
    \end{tabular}    \caption{Percentage of accidents divided by the range of losses in load values.}
    \label{tab:table_5}
\end{table}

\begin{figure}
    \centering
    \includegraphics[width=0.9\textwidth]{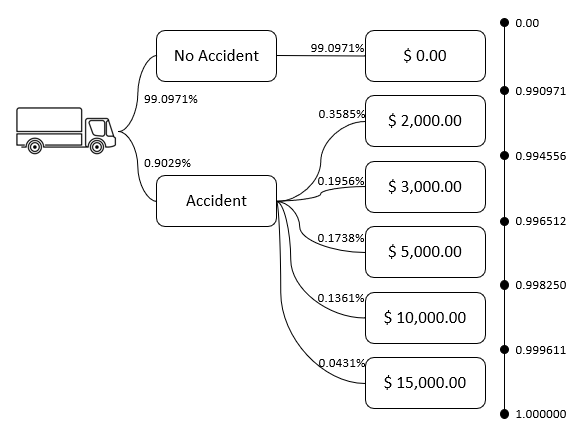}
    \caption{Accident probabilities and its costs to Monte Carlo Simulation.}
    \label{figure6}
\end{figure}

The Figure \ref{figure6} illustrates an example of how the probabilities were distributed. The maximum load value for each range and the cost to the road freight company was considered as being 1\% of its value, a deductible that load insurances usually charge. Thus, the occurrence of the event in the range between \$0.01 and R\$200,000.00, the value to be considered will always be the highest of the range, which in this case will be \$200,000.00 and deductible cost for the carrier would be \$2,000.00.

Also through Figure \ref{figure6} it is possible to verify the accumulated percentage values, between 0 and 1, for each accident cost on the right. This will be important for the Monte Carlo simulation, which at each iteration selects a random value between 0 and 1 that will correspond to a accident cost. For example, according to Figure \ref{figure6} each percentage range is equivalent to its cost, so any value selected in the range between 0 and 0.990971 will correspond to a accident cost equal to \$0.00.

In this way, 1,000,000 iterations were performed for each arc of the problem and the average of accident costs was calculated and it generates $r_{ij}$ that will be used in the objective function of this VRP. This number of iterations was chosen because the accidents probabilities are very low and thus the values of $r_{ij}$ present better convergence. However, if the iterations are above the established, the resolution time becomes longer and therefore, for this problem, 1,000,000 is an ideal number.

\section{Results and Discussions}
\subsection{Probabilities and costs of risks}

\begin{figure}
    \centering
    \includegraphics[width=1\textwidth]{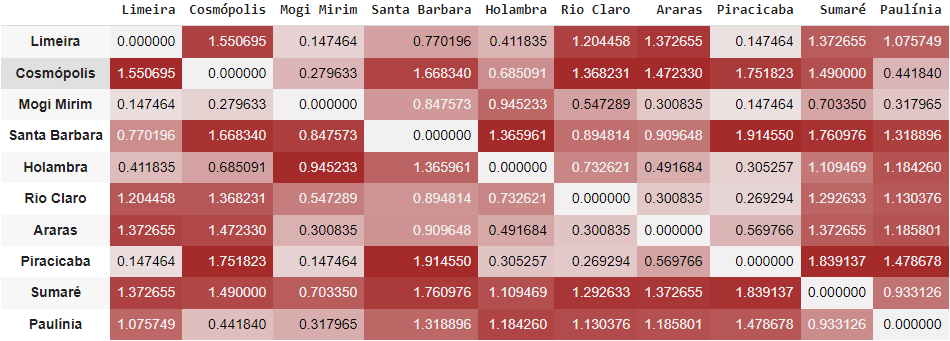}
    \caption{Accidents Probabilities $Paccident_{ij}$ in \%.}
    \label{figure7}
\end{figure}

\begin{figure}
    \centering
    \includegraphics[width=1\textwidth]{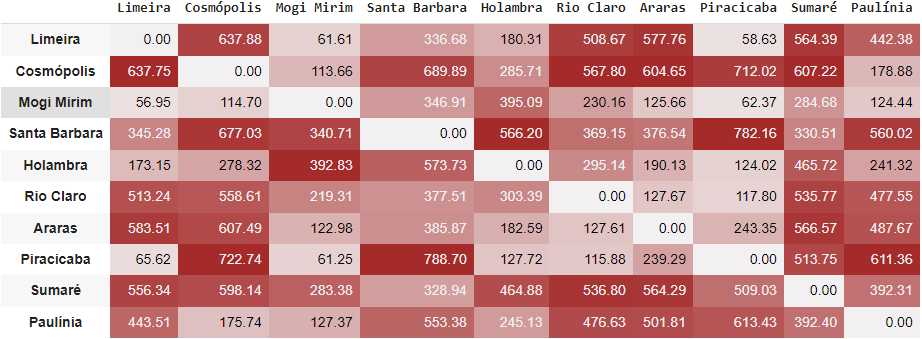}
    \caption{Risk cost $r_{ij}$ in \$.}
    \label{figure8}
\end{figure}

Figures \ref{figure7} and \ref{figure8} represent the results for $Paccident_{ij}$ and $r_{ij}$, respectively. After performing analyses, it was observed that $Paccident_{ij}$ and $r_{ij}$ obtained are in line with what was expected. When comparing the arcs of the same type of road as Piracicaba-Santa Bárbara and Limeira-Mogi Mirim, it was verified that $Paccident_{ij}$ for the first (1.91455 \%) is greater than the second (0.147464 \%) as well as for $r_{ij}$ which is \$ 788.70 and \$ 61.61, respectively. It was expected because the vehicles flow on the Piracicaba-Santa Barbara arc (6,943 heavy vehicles) is greater than the other (535 heavy vehicles). 

When comparing arcs that present similar or close vehicle flow such as Mogi Mirim-Araras and Mogi Mirim-Limeira, it is noted that $Paccident_{ij}$ for the first (0.300835 \%) is greater than the second (0, 147464 \%) as well as for $r_{ij}$ which is \$ 125.66 and \$ 56.95. And this is due to the fact that the road that connects Mogi Mirim and Araras is a single lane two way road, according to the nomenclature of \cite{cnt_painelacidentes2018}, which has the highest death rate among all type of roads.

Other analyzes like this were also carried out and the same conclusions were reached. In this way, it is possible to state that the methodology followed to find $Paccident_{ij}$ and $r_{ij}$ was satisfactory, since consistent results were found according to the results obtained.

\subsection{Optimization results}

\begin{figure}
    \centering
    \includegraphics[width=1.0\textwidth]{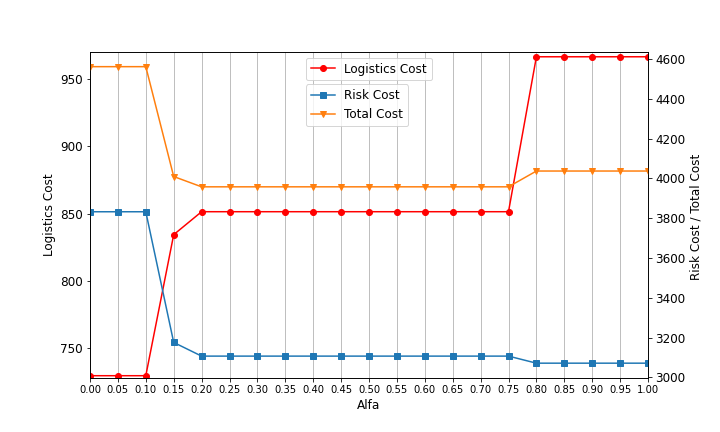}
    \caption{Results of $c_{ij}$, $r_{ij}$ e $z$}
    \label{figure9}
\end{figure}

Figure \ref{figure9} represents the graph with the results of each parcel, logistics and risk cost, and the sum between them for the optimization varying $\alpha$. Note that for $\alpha = 0$, the logistic cost assumes the minimum value when compared to others $\alpha$, but the risk is the highest. When considering a greater weight of security during the optimization, that is, when increasing $\alpha$, the risk cost starts to decrease and the logistic cost increases. This occurs because the model starts to prioritize arcs with lower risks and thus reducing the concern with the cost of tolls and fuel.

It is verified by Figure \ref{figure9} that for this problem a variation is perceived in $\alpha = 0.15$ which the logistic cost suffers an increase from \$ 729.70 to \$834.36 while the risk decreases from \$3832.44 to \$3175.29 representing a decrease of approximately 17.15 \%. Between $\alpha = 0.20$ and $\alpha = 0.75$, the minimum logistics costs remain at \$851.38 and risks at \$3106.47, which further increases the level of safety fact.

It is important to a decision maker to see it because the choice of routes will depend on what level of security it wants to work for its fleet. If security is considerable in the process of the company, it is important to work with $\alpha \geq 0.15$ but if the logistics cost is still relevant, the level of security should be $\alpha \leq 0.75$, because for $\alpha \geq 0.80$ the logistics cost will reach R\$ 966.34 and the risk R\$ 3070.95 which allows us to say that, due to the first parcel, there will be a considerable increase of approximately 13.50 \% and the second reduce only 1.1 \%, it would not be worth working in this range.

\begin{figure}
    \centering
    \subfigure[subfigure1][$\alpha = 0,00$]{\includegraphics[width=0.47\textwidth]{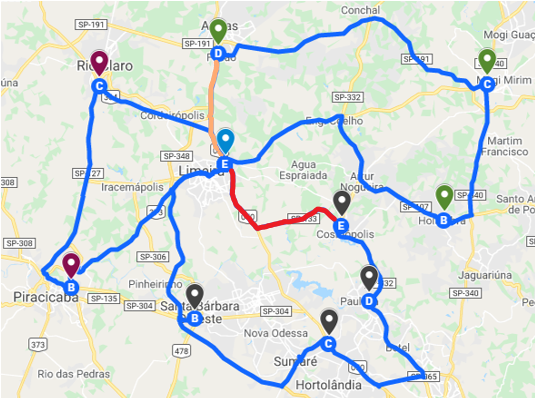} \label{subfigure1}}
    \hfill
    \subfigure[subfigure2][$\alpha = 0,15$]{\includegraphics[width=0.47\textwidth]{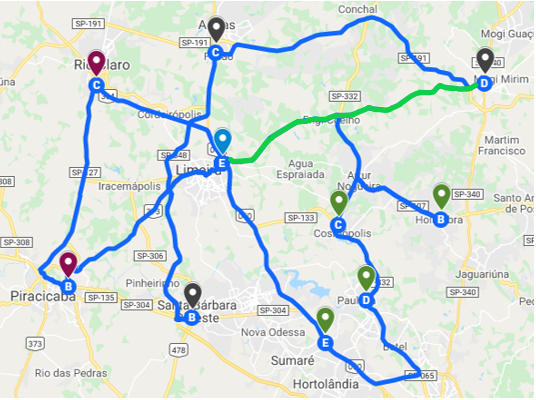} \label{subfigure2}}
    \hfill
    \subfigure[subfigure3][$\alpha = 0,20$]{\includegraphics[width=0.47\textwidth]{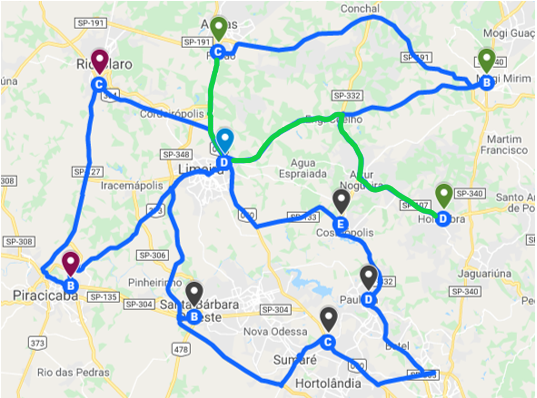}\label{subfigure3}}
    \hfill
    \subfigure[subfigure4][$\alpha = 1,00$]{\includegraphics[width=0.47\textwidth]{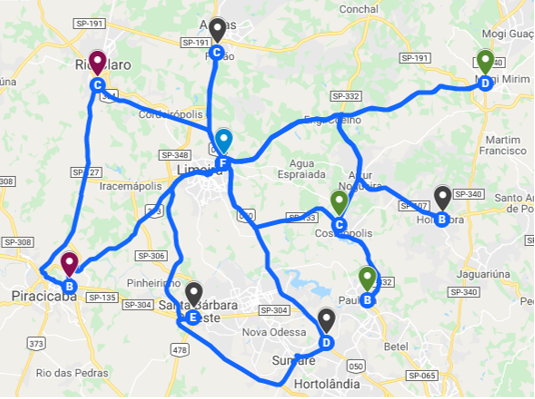}\label{subfigure4}} \hfill    
    \caption{Routes obtained for different values of $\alpha$.}
    \label{figure10}
\end{figure}

The routes were also evaluated and it is represented by Figures \ref{subfigure1}, \ref{subfigure2}, \ref{subfigure3} and \ref{subfigure4} which each delivery point is identified with a marker and letter, referring to the sequence of deliveries made by each truck. It is also noted that there are different colors for the markers that identify which truck is carrying out the delivery, facilitating the evaluation of the results. For example, the markers in red represent truck 1, in black truck 2, in green truck 3 and in blue the warehouse.

Only routes were printed from $\alpha$ being 0.00, 0.15, 0.20 and 1.00. According to the graph in Figure \ref{figure9}, these are the points that show variations. For $0.10 \leq \alpha \leq 0.15$ there were some variations, the main one being eliminating the Cosmópolis-Limeira ($r_{ij} = 637.75$ red arc in Figure \ref{subfigure1}) and adding Mogi Mirim-Limeira ( $r_{ij} = 61.61$ green arc in the Figure \ref{subfigure2}) and thus reducing the risk, in this case alone, by approximately 85\% in relation to the total amount that was reduced for this comparison.

The arc Mogi Mirim-Limeira or vice versa, represented by green arc in Figure \ref{subfigure2}, it started to be added as soon as $\alpha$ increases and it occurs because its risk cost is one of the lowest. Its inclusion is the one that most impacted the minimization of risk, therefore for $\alpha \geq 0.20$ the risk cost decreases only approximately to 1\%, as this arc has already been considered.

When analyzing $\alpha = 0.20$, it was noticed that the arc Araras-Holambra (green arc in the Figure \ref{subfigure3}) was included ($r_{ij} = 182.59$) but a part of it passes through the same road with origin in Araras and destination Limeira ($ r_{ij} = 583.51$), represented by orange arc in Figure \ref{subfigure1}. This means that to leave Araras and go to Holambra, the truck initially travels along a less safe road, but then in the most part of the arc is on low risk roads, which makes $r_{ij}$ to have a lower value and this could be a limitation for the method.

\section{Conclusion}

This paper consisted of inserting the issue of route safety in the vehicle routing process of a road freight company in order to help the decision maker to select the best routes that minimize both the distance and the risk of accidents.

A method was developed that considers the statistical analysis to estimate the accidents probabilities for each arc and the costs related to the risks of accidents from the Monte Carlo simulation. The results obtained were coherent, converging to what was expected to the problem.

Through a mathematical model of VRP, logistics costs and risks were minimized and the results were analysed for each part of the objective function and the routes obtained for different values of a risk coefficient alfa ($\alpha$).

According to the outcomes displayed in Figure \ref{figure9}, it was identified that as the safety factor increases, the risk-related costs decrease, which is more noticeable when $0.10 \leq \alpha \leq 0.20$. It was also noticed that, when comparing the routes, it was identified that the model started to incorporate security routes as results.

A limitation observed in the proposed method is for two trips passing through the same road, the accidents probabilities were different and this can be explained because in one trip the truck would be on the more dangerous road for longer than the other, interfering estimation process.

However, it is possible to analyze that the method worked in a positive way by verifying the accidents probabilities and risk cost and it can help the decision maker of a road freight company to select the best routes considering distance and accident risk. In addition, the whole approach developed here is simple and adaptable to any VRP model and can be used by any company for free. 

Finally \textit{Knime Analytics Platform} helps to deal with the use of real data in this paper simplifying data exploration, analysis, visualization and interpretation.

For future works, it is relevant to develop an approach that considers the risk of load theft in VRP and based on statistical analysis. In addition, it is necessary to add that during the day accident probability and risk of theft could vary because some variables, as traffic flow, changes according the time and this will help road freights and load insurers to make decisions about the time that a truck could pass at certain road.

%
%

\bibliographystyle{splncs04}

\end{document}